\documentclass{amsart}

\newtheorem{theorem}{Theorem}[section]

\newtheorem{lemma}[theorem]{Lemma}

\newtheorem{conjecture}[theorem]{Conjecture}

\theoremstyle{definition}
\newtheorem{definition}[theorem]{Definition}
\newtheorem{propdef}[theorem]{Proposition-Definition}

\theoremstyle{remark}
\newtheorem{remark}[theorem]{Remark}

\numberwithin{equation}{section}

\newcommand{\PL}{{\mathbb{P}^1}}
\newcommand{\AL}{{\mathbb{A}^1}}
\newcommand{\CT}{{\mathbb{C}^\times}}
\newcommand{\PP}{{\mathbb{P}^2}}

\begin{document}

\title{Log mirror symmetry and local mirror symmetry}
\author{Nobuyoshi Takahashi}
\address{Department of Mathematics, 
Hiroshima University, Higashi-Hiroshima 739-8526, Japan}
\email{takahasi@math.sci.hiroshima-u.ac.jp}

\begin{abstract}
We study Mirror Symmetry of log Calabi-Yau surfaces. 
On one hand, 
we consider the number of ``affine lines'' of each degree 
in $\PP\backslash B$, where $B$ is a smooth cubic. 
On the other hand, 
we consider coefficients of a certain expansion of a function 
obtained from the integrals of $dx/x\wedge dy/y$ 
over 2-chains whose boundaries lie on $B_{\phi}$, 
where $\{B_{\phi}\}$ is a family of smooth cubics. 
Then, for small degrees, they coincide. 

We discuss the relation between this phenomenon 
and local mirror symmetry for $\PP$ in a Calabi-Yau 3-fold
(\cite{CKYZ}). 
\end{abstract}

\maketitle

In the classification theory of algebraic varieties, 
one can study non-compact algebraic manifolds 
by means of Log Geometry. 

For a quasi-projective manifold $U$, 
there exist a projective variety $X$ 
and a normal crossing divisor $B$ on $X$ 
such that $X\backslash B$ is isomorphic to $U$. 
A rational differential form on $X$ is 
called a logarithmic form 
if it can be locally written as a linear combination 
(with regular functions as coefficients) 
of products of $d\log f$, 
where $f$ is a regular function that vanishes at most in $B$. 
Then the space of logarithmic forms is an invariant of $U$ 
and it plays the role that the space of regular forms plays 
in the classification theory of projective manifolds. 

In particular, we may think of a pair $(X, B)$ 
(or the open variety $X\backslash B$) 
as a log Calabi-Yau manifold 
if $K_X+B$ is trivial. 
It is expected that some analogue of Mirror Symmetry holds 
for log Calabi-Yau manifolds. 

As an evidence, 
we studied Mirror Symmetry 
for the one-dimensional log Calabi-Yau manifold, 
i.e. $(\PL, \{0,\infty\})$ or $\CT$, in \cite{T2}: 
for $b\geq0$ and $k, l>0$, 
we consider covers of $\PL$ 
that have $k$ and $l$ points over $0$ and $\infty$ 
and are simply branched over $b$ prescribed points. 
We define $F_{b,k,l}(z_1, \dots, z_k; w_1, \dots, w_l)$ 
to be the generating function of the number of such curves. 
Then $F_{b,k,l}$ coincides with the sum of certain integrals 
associated to graphs 
that have $k$ and $l$ `initial' and `final' vertices, 
have $b$ internal vertices 
and are trivalent at internal vertices 
(here, variables $z_1, \dots, z_k$ and $w_1, \dots w_l$ 
are attached to initial and final vertices). 

The aim of this paper is to study the two-dimensional case 
as a further example of Log Mirror Symmetry. 

In Section \ref{log}, 
we describe Mirror Symmetry of $\PP\backslash B$ 
where $B$ is a smooth cubic. 
On A-model side, 
we count curves of each degree in $\PP\backslash B$ 
whose normalizations are $\AL$ 
and construct a power series from those numbers. 
On B-model side, we construct a function 
from the integrals of $dx/x\wedge dy/y$ 
over 2-chains whose boundaries lie on $B_{\phi}$, 
where $B_{\phi}$ is a member 
of a family of smooth cubics parametrized by $\phi$. 
Then we see that the coefficients of their expansions coincide 
up to order 8. 

Then, in Section \ref{log_and_local}, 
we discuss the relation between our log mirror symmetry 
and local mirror symmetry studied in \cite{CKYZ}. 
We explain that the number of ``affine lines'' 
is in some sense the dual of 
local Gromov-Witten invariants of 
$\PP$ in a Calabi-Yau 3-fold.

\section{Log mirror symmetry}\label{log}

\subsection{A-model}
On the enumerative side, 
we counted curves of the lowest `log genus' in \cite{T1}. 

\begin{definition}
A plane curve $C$ is said to satisfy the condition (AL) 
if it is irreducible and reduced 
and the normalization of $C\backslash B$ 
is isomorphic to the affine line $\AL$. 
\end{definition}

\begin{remark}
If $C$ satisfies (AL), 
then $C\cap B$ consists of one point. 
This point is a $(3\deg C)$-torsion 
for a group structure on $B$ 
whose zero element is an inflection point of $B$. 
\end{remark}

We treat only `primitive' cases. 

\begin{definition}
Let $P$ be a point of \emph{order} $3d$ 
for a group structure on $B$ whose zero element is a point of inflection. 

Then we define $m_d$ to be the number of curves $C$ of degree $d$ 
which satisfy (AL) and $C\cap B=\{P\}$
(or, more precisely, the degree of the 0-dimensional scheme 
parametrizing such curves). 
\end{definition}

The above definition implicitly assumes 
that the number $m_d$ does not depend on the choice of $P$. 
Although we haven't proved it yet, 
it holds in the cases where we know $m_d$, i.e. when $d\leq8$. 

\begin{theorem}(\cite{T1})\label{table}
We have the following table for $m_d$: 
\[
\begin{tabular}{|c||c|c|c|c|c|c|}
  \hline
  $d$   & $1$ & $2$ & $3$ &  $4$ &   $5$ &   $6$ \\
  \hline
  $m_d$ & $1$ & $1$ & $3$ & $16$ & $113$ & $948$ \cr
  \hline
\end{tabular}
\]
Furthermore, under a technical hypothesis(see \cite{T1}), 
we have $m_7=8974$ and $m_8=92840$. 
\end{theorem}

Our invariants are apparently related 
to the relative Gromov-Witten invariants 
defined in \cite{IP} and \cite{LR} for pairs of a symplectic variety 
and its subvariety of real codimension 2. 
In algebraic language: 

\begin{propdef}\label{relgw}(1-pointed case of \cite{G})
Let $\bar{M}_{0,1}(\PP,d)$ be the moduli stack of 
$1$-pointed genus 0 stable maps of degree $d$ to $\PP$ 
and $\bar{M}_i^B(\PP,d)$ the closed subset 
consisting of points corresponding to $f:(C,P)\to\PP$ 
such that 

(i) $f(P)\in B$ and 

(ii) $f^*B-iP\in A_0(f^{-1}B)$ is effective. 

Then the virtual fundamental class $[\bar{M}_i^B(\PP,d)]^{virt}$ 
is naturally defined and is of expected dimension $3d-i$. 
\end{propdef}

Considering that the number of $3d$-torsions is $(3d)^2$ 
and that the relative Gromov-Witten invariants 
take multiple covers into account, 
we expect the following to hold. 

\begin{conjecture}\label{conjrelgw}
$[M_{3d}^B(\PP,d)]^{virt}=(3d)^2\sum_{k|d}(-1)^{d-d/k} m_{d/k}/k^4$. 
\end{conjecture}

\begin{remark}
(1) Although the factor $k^{-4}$ looks unfamiliar, 
this is compatible with the factor $k^{-3}$ 
for Gromov-Witten invariants, 
since $m_d$ is conjecturally equal to $(-1)^{d-1}n_d/(3d)$, 
where $n_d$ is the local Gromov-Witten invariant
(see Remark \ref{remlogandlocal}). 

(2) A.~Gathmann informed the author that this is true 
for $d\leq8$(i.e. where $m_d$ is calculated). 
\end{remark}

\subsection{B-model}
The result of \cite{T2} suggests 
that the mirror manifold of $\CT$ is $\CT$. 
In A-model, we associate a parameter 
to each point over $0$ or $\infty$ on a cover of $\CT$, 
and in B-model this corresponds to the choice of points on $\CT$: 
the former can be considered as ``K\"ahler moduli'' 
and the latter as ``complex moduli''. 
In B-model, we used
the points as the boundary of integrals 
to define coordinate functions. 

With this in mind, 
we look for a function which has an expansion with coefficients $m_d$. 
Although the calculation below is 
essentially the same as in \cite{CKYZ}, 
here we see it from the point of view of Mirror Symmetry 
of log Calabi-Yau surfaces: 
we claim that the mirror of 
($\PP\backslash(\hbox{smooth cubic})$ \& K\"ahler moduli) 
is ($\CT\times\CT$ \& cubic), 
where the latter cubic has complex moduli 
and we take it as the boundary of integrals, 
just as in the case of $\CT$. 

We consider homogeneous coordinates $X, Y, Z$ 
and inhomogeneous coordinates $x, y$ on $\PP$. 
Let $\Omega$ be the logarithmic 2-form $dx/x\wedge dy/y$ 
on $X_0:=\PP\backslash\{XYZ=0\}$ 
and write $B_\phi$ for the plane cubic 
defined by $XYZ-\phi(X^3+Y^3+Z^3)=0$. 

Then we consider the integral 
\[
  I:=\int_{\Gamma_\phi} \Omega, 
\]
where $\Gamma_\phi$ is a 2-chain in $X_0$ 
whose boundary has support on $B_\phi$. 

\begin{lemma}
We have 
\[
  \phi\frac{dI}{d\phi} = \int_{\partial\Gamma} dx/(x-3\phi y^2). 
\] 
\end{lemma}
\begin{proof}
If we fix $x$ and differentiate $xy-\phi(x^3+y^3+1)=0$, 
we have $dy/d\phi=xy/\phi(x-3\phi y^2)$. 
\end{proof}

So, if we set $z:=\phi^3$ and $\theta:=z\frac{d}{dz}$, we have 
\[
  \{\theta^3-3z\theta(3\theta+1)(3\theta+2)\}I=0. 
\]
The following functions form a basis of the space of solutions: 
\begin{eqnarray*}
I_1 & = & 1, \\
I_2 & = & \log z + I^{(0)}_2, \\
I_3 & = & I_2 \log z -\frac{(\log z)^2}{2} + I^{(0)}_3, 
\end{eqnarray*}
where 
\begin{eqnarray*}
I^{(0)}_2 
  & = & 6 z + 45 z^2 + 560 z^3 + \frac{17325}{2} z^4 
  + \frac{756756}{5} z^5 + 2858856 z^6 \\ 
  &  & + \frac{399072960}{7} z^7 + \frac{4732755885}{4} z^8 + \cdots, \\
I^{(0)}_3 
  & = & 9 z + \frac{423}{4} z^2 + 1486 z^3 + \frac{389415}{16} z^4 
  + \frac{21981393}{50} z^5 + \frac{16973929}{2} z^6 \\ 
  &   & + \frac{8421450228}{49} z^7 + \frac{1616340007953}{448} z^8 
  + \cdots. 
\end{eqnarray*}

The monodromy of $I_3$ for $z\to e^{2\pi i}z$ is 
$2\pi i(I_2-\pi i)$, which we denote by $\tilde{I}_2$. 
We denote the monodromy $(2\pi i)^2$ of $\tilde{I}_2$ by $\tilde{I}_1$. 
Now we write $I_3$ in terms of 
$q:=e^{2\pi i\tilde{I}_2/\tilde{I}_1}=-e^{I_2}=-ze^{I^{(0)}_2}$: 
\begin{eqnarray*}
I_3 
  & = & \frac{(I_2)^2}{2} + 9 z + \frac{351}{4} z^2 + 1216 z^3  
   + \frac{319455}{16} z^4 + \frac{18122643}{50} z^5  \\
  &   &  + \frac{35161224}{5} z^6 + \frac{7009518168}{49} z^7 
   + \frac{1350681750297}{448} z^8 + \cdots \\
  & = & \frac{(\log(-q))^2}{2}  + 0 - 9 q + \frac{135}{4} q^2 
   - 244 q^3 + \frac{36999}{16} q^4 - \frac{635634}{25} q^5 \\
  &   &  + 307095 q^6 - \frac{193919175}{49} q^7 
   + \frac{3422490759}{64} q^8 + \cdots \\
      & = & \frac{(\log(-q))^2}{2} - 3^2.1 \sum_{k=1}^\infty \frac{q^k}{k^2} 
            + 6^2.1 \sum_{k=1}^\infty \frac{q^{2k}}{k^2} 
            - 9^2.3 \sum_{k=1}^\infty \frac{q^{3k}}{k^2} 
            + 12^2.16 \sum_{k=1}^\infty \frac{q^{4k}}{k^2} \\
      &   & - 15^2.113 \sum_{k=1}^\infty \frac{q^{5k}}{k^2} 
            + 18^2.948 \sum_{k=1}^\infty \frac{q^{6k}}{k^2} 
            - 21^2.8974 \sum_{k=1}^\infty \frac{q^{7k}}{k^2} \\
      &   & + 24^2.92840 \sum_{k=1}^\infty \frac{q^{8k}}{k^2} + \cdots. 
\end{eqnarray*}

\begin{remark}
$K:=(zdI_2/dz)^3 d^2I_3/dI_2^2$ 
satisfies $dK/dz=(27/(1-27z))K$, 
and therefore we have $K=1/(1-27z)$. 
\end{remark}

Comparing the coefficients with the table in Theorem \ref{table}, 
we propose: 

\begin{conjecture}\label{conjmirror}
\[
  I_3 = \frac{(\log(-q))^2}{2} + 
    \sum_{d=1}^\infty (-1)^d(3d)^2 m_d \sum_{k=1}^\infty \frac{q^{dk}}{k^2}. 
\]
\end{conjecture}

\begin{remark}
If we assume Conjecture \ref{conjrelgw}, 
the previous conjecture is equivalent to 
\[
  I_3 = \frac{(\log(-q))^2}{2} + 
    \sum_{d=1}^\infty (-1)^d [M_{3d}^B(\PP,d)]^{virt} q^d,  
\]
and in fact this may be a more natural equality. 

According to A.~Gathmann, his algorithm in \cite{G} 
can be used to prove this. 
\end{remark}

\section{Log Mirror and Local Mirror}\label{log_and_local}

In \cite{CKYZ}, the generating function of 
the ``numbers of rational curves in a local Calabi-Yau 3-fold'' 
was given. 

Let $\PP$ be embedded in a Calabi-Yau 3-fold
and denote by $n_d$ the contribution of rational curves of degree $d$ 
in $\PP$ to the number of rational curves in the Calabi-Yau 3-fold. 
Let $\bar{M}_{0,0}$ be the moduli of stable maps 
of genus 0 curves to $\PP$ with degree $d$ images 
and $U$ the vector bundle over $\bar{M}_{0,0}$ 
whose fiber at the point $[f:C\to\PP]$ 
is $H^1(C, f^*K_\PP)$.  
Then, the Chern number $K_d:=c_{3d-1}(U)$ 
is equal to $\sum_{k|d}n_{d/k}/k^3$.

\begin{theorem}(\cite{CKYZ})
\begin{eqnarray*}
  I_3 
  & = & 
    \frac{(\log(-q))^2}{2} - 
    \sum_{d=1}^\infty 3d K_d q^d \\
  & = &
    \frac{(\log(-q))^2}{2} - 
    \sum_{d=1}^\infty 3d n_d \sum_{k=1}^\infty \frac{q^{dk}}{k^2}. 
\end{eqnarray*}
\end{theorem}

\begin{remark}\label{remlogandlocal}
Thus Conjecture \ref{conjmirror} 
is equivalent to $n_d = (-1)^{d-1}3dm_d$, 
and this equality holds for $d\leq8$ 
by Theorem \ref{table} 
and the $q$-expansion of $I_3$ in the previous section. 
\end{remark}

We give a heuristic argument as to why this should hold, 
although it contains serious gaps as explained later. 

By Serre duality, the dual of $U$ is isomorphic to 
the vector bundle $V$ over $\bar{M}_{0,0}$ 
whose fiber at $[f:C\to\PP]$ is 
$H^0(C, K_C\otimes f^*\mathcal{O}_\PP(B))$.  
Since the rank of the bundles is $3d-1$, 
we have $c_{3d-1}(V)=(-1)^{d-1}K_d$. 

Let $\bar{M}_{0,1}$ be the moduli of stable maps 
of 1-pointed genus 0 curves to $\PP$ with degree $d$ images, 
$M_{0,1}$ the open subset of $\bar{M}_{0,1}$ 
representing $f:(\PL,P)\to\PP$ 
and $\pi:\bar{M}_{0,1}\to \bar{M}_{0,0}$ the projection. 
Further, let $E_1$ be the line bundle over $\bar{M}_{0,1}$ 
whose fiber at $f:(C,P)\to\PP$ is 
$H^0(C, f^*\mathcal{O}_\PP(B)\otimes\mathcal{O}_P)$ 
and $L$ the line bundle whose fiber is 
$H^0(C, \mathcal{O}_P(-P))$. 
Then, we have $c_{3d}(\pi^*V\oplus E_1) = 3dc_{3d-1}(V)$, 
for the zero set of a section of $E_1$ 
induced by a defining equation of $B$ 
is the set of the points corresponding to $f:(C,P)\to\PP$ 
such that $f(P)\in B$, 
and there are $3d$ such points for any $f:C\to\PP$. 

We have an exact sequence 
\[
  0\to K_C\to K_C(P)\overset{\hbox{\tiny residue}}\to\mathcal{O}_C\to0. 
\]
If $C$ is irreducible, i.e. isomorphic to $\PL$, 
we obtain exact sequences 
\begin{eqnarray*}
  0 & \to & H^0(C,K_C(-(k+1)P)\otimes f^*\mathcal{O}_\PP(B)) 
            \to H^0(C,K_C(-kP)\otimes f^*\mathcal{O}_\PP(B)) \\
    & \to & H^0(C,O_P(-(k+1)P)\otimes f^*\mathcal{O}_\PP(B)) \to 0
\end{eqnarray*}
for $0\leq k\leq3d-2$. 
We also have $H^0(C,K_C(-(3d-1)P)\otimes f^*\mathcal{O}_\PP(B))=0$. 

Thus, on $M_{0,1}$, 
we have a filtration of $\pi^*V\oplus E_1$ 
such that the associated graded module is isomorphic to 
$\bigoplus_{i=0}^{3d-1}E_1\otimes L^i$. 

On the other hand, 
there are about $(3d)^2m_d$ plane curves 
of degree $d$ satisfying (AL), 
since the number of $3d$-torsions on $B$ is $(3d)^2$. 
They are in one-to-one correspondence 
with points $(f:(\PL,P)\to\PP)\in M_{0,1}$ 
such that $f$ is birational onto the image 
and that $f^*B=3dP$. 

Consider the vector bundle $E_{3d}$ of rank $3d$ 
over $\bar{M}_{0,1}$ 
whose fiber at $f:(C,P)\to\PP$ is 
$H^0(C, f^*\mathcal{O}_\PP(B)\otimes\mathcal{O}_{3dP})$. 
On $M_{0,1}$, 
the zero set of a section of $E_{3d}$ 
induced by a defining equation of $B$ 
is the set of the points $f:(C,P)\to\PP$ 
such that $f^*B=3dP$. 

Now, from the exact sequences 
\[
  0\to\mathcal{O}_P(-kP)\to\mathcal{O}_{(k+1)P}\to\mathcal{O}_{kP}\to0, 
\]
we see that 
$E_{3d}$ has a filtration 
such that the associated graded module is isomorphic to 
$\bigoplus_{i=0}^{3d-1}E_1\otimes L^i$. 
Thus we may expect 
$(3d)^2m_d\approx(-1)^{3d-1}3dK_d\approx(-1)^{3d-1}3dn_d$. 

\bigskip
Rigorously, however, this argument makes little sense. 
First, the section of $E_{3d}$ 
induced by a defining equation of $B$ 
has undesirable zeros 
in $\bar{M}_{0,1}\backslash M_{0,1}$. 
For example, if $C=C_1\cup C_2\cup C_3$(a chain in this order), 
$P\in C_2$ and $C_2$ maps to a point $Q$ in $B$, 
we may take any rational curves through $Q$ 
as the images of $C_1$ and $C_3$. 
Thus the number $(3d)^2m_d$ may be much different 
from $c_{3d}(E_{3d})$. 
Second, we have a filtration of $\pi^*V\oplus E_1$ 
merely on $M_{0,1}$. 

Let $\bar{M}_i^B(\PP,d)$ be as in Definition \ref{relgw}. 
The section of the line bundle $(E_1\otimes L^i)|_{\bar{M}_i^B(\PP,d)}$ 
induced by a defining equation of $B$ 
vanishes when $f^{-1}B\supseteq (i+1)P$ is satisfied. 
Then, \cite{G} describes the difference between 
$c_1(E_1\otimes L^i).[\bar{M}_i^B(\PP,d)]^{virt}$ 
and $[\bar{M}_{i+1}^B(\PP,d)]^{virt}$.
This should account for the difference 
between $c(\bigoplus_{i=0}^{3d-1}E_1\otimes L^i)$ and 
$c(\pi^*V\oplus E_1)$.

\end{document}